\documentclass[preprint,12pt]{elsarticle}

\usepackage{amssymb}
\usepackage{bm}
\usepackage{amsmath}
\usepackage{cleveref}
\usepackage{extarrows}
\usepackage{subfigure}
\usepackage{color} 
\usepackage{boxedminipage}
\usepackage{diagbox} 
\usepackage{multirow}
\usepackage[sc]{mathpazo}
\newcommand{\tabincell}[2]{\begin{tabular}{@{}#1@{}}#2\end{tabular}}
\usepackage{nomencl}
\makenomenclature

\journal{}

\begin{document}

\begin{frontmatter}



\title{Efficient calculation of the integral equation for simulating 2D TE scattering in a homogeneous medium using the Ewald method and a Gabor frame discretization}


\author[label1]{Xinyang Lu}
      \address[label1]{School of Mathematical Sciences, Zhejiang University, Hangzhou 310027, Zhejiang, China}
       \address[label2]{Department of Electrical Engineering, Eindhoven University of Technology, P.O. Box 513, 5600 MB Eindhoven, The Netherlands}
\author[label2]{M. C. van Beurden}
\author[label1]{Qingbiao Wu}

\begin{abstract}
We utilize the domain integral equation formulation to simulate two-dimensional transverse electric scattering in a homogeneous medium and a summation of modulated Gaussian functions to approximate the dual Gabor window.
Then we apply Ewald Green function transformation to separate the integrals related to $x$ and $z$ in the integral equation, which 
produce Gaussian functions. These Gaussian functions in the integrands can be integrated analytically, which greatly simplifies the calculation process. Finally, we discuss the convergence and the selection of the Ewald splitting parameter $\mathcal{E}$.
\end{abstract}


\end{frontmatter}


\printnomenclature[1in]

\section{Introduction}
\label{introduction}
 \text{The} domain integral equation approach {has been used} to solve the electromagnetic scattering problems in recent decades \cite{PP1991Aweak, AIM1996EB, Lalanne2007Numerical}. {The Conjugate Gradient Fast Fourier Transform (CGFFT) \cite{Zwamborn1992The} and pre-corrected FFT \cite{Phillips1994Proceedings} have been efficient ways to solve the integral equation.}
 The Green function play{s} an important role in this method. How to overcome the singularity and coupling of the variables in the Green function is the key to solve the integral equation. Dilz and {Van Beurden} \cite{Dilz2016The} proposed 
the mixed spatial and spectral approach with Gabor frame discretization. It utilizes the Fourier transformation to eliminate the coupling and deal with the singularity in spectral domain by scaling coordinates. The numerical results demonstrate that it
yields an efficient and accurate way to numerically simulate the scattering from a finite dielectric object in a homogeneous background medium. A Gabor frame \cite{Bastiaans1981A, Bastiaans1995Gabor} is rapidly convergent both in spatial and spectral domain, which makes it possible for us to use the mixed spatial and spectral method to solve the domain integral function efficiently. 

The Gabor frame is {presented} in \cite{Feichtinger1998Gabor} and there is a brief introduction in Section 2.2, in which the formulas show that the dual window $\eta(x)$ plays an important role in the calculation of Gabor coefficients. Since the Gabor coefficients are not uniquely determined because of oversampling, the choice for the dual window function is not unique. Zak transform with generalized Moore-Penrose pseudo inverse is {a} popular way to calculate $\eta(x)$ \cite{Bastiaans1995Gabor, Werther2005Dual} and $\eta$ obtained via this method  is the optimum solution in the sense of minimum $\mathcal{L}^2$ norm \cite{Feichtinger1998Gabor, Janssen1994Signal}. However we can only obtain discrete values of $\eta(x)$, which makes it difficult to derive analytical results of these integrals related to $\eta$.  Here, we assume that 
{the dual window is represented via the original Gabor frame of the expansion.} Then we can simplify the integral containing $\eta$ by using the property of Gaussian functions.

In the integral equation, we must deal with the singularity problem of {the} Green's function.  Although a coordinate scaling to smoothen the branchcut of the Green function in the spectral domain works well \cite{Dilz2016The}, we want {a} faster and more accurate means for evaluating the Green function in the spatial domain. Ewald Green's function transformation, proposed by P. P. Ewald and described in \cite{Ewald1921Die, Ewald1970On},  provides the integral formula of the {free-space} Green's function for {three-dimensional} (3D) problems \cite{Jordan1986An} and 2D problems \cite{Mathis1996A}. Applying this formula to the integral equation, the integrals related to $x$ and $z$ are separated, which greatly reduces the difficulty of the integral. Meanwhile, the integral related to $x$ also has the Gaussian integral form and it can be combined with $\eta$ to get some analytic results. 

To {handle} the singularity of Green function, we use the Ewald method \cite{Capolino2005Efficient, Capolino2007Efficient,  Komanduri20161} that splits the integral {representation} of Green function into two parts. The selection of the Ewald splitting parameter $\mathcal{E}$ {is} discussed. 
Only $G_{spectral}$ suffers from the singularity in spatial domain and we utilize the Fourier {transformation} to solve $G_{spectral}$ in spectral domain. Finally, what we need to calculate is several {one-dimensional} integrals for the products of the continuous functions. The numerical examples indicate that the results obtained from our proposed method match well that in \cite{Dilz2016The} with less computation and storage.
\section{Statement of the problem}
\label{Statement of the problem}
The derivation of the problem formulation can be found in  \cite{Dilz2016The} and \cite{RJDilz2017spatialspectral}. We will only present the main formulas used in the calculation and more details can be found in Chapter 2 and Chapter 5 in \cite{RJDilz2017spatialspectral}. 

\subsection{Problem formulation}
For two dimensional transverse electric scattering in a homogeneous isotropic dielectric medium, the final equation we will solve is as follows:
 \begin{small}
 \begin{equation}
 \begin{split}
 &\qquad \chi(x,z)E^s(x,z) = \\
 & k^2_0\varepsilon_{rb}\chi(x,z)\int \limits_x\int \limits_z G(x,z|x',z')\cdot\chi(x',z') \cdot  [E^i(x',z') + E^s(x',z')]dx'dz',
 \end{split}
 \end{equation}
 \end{small}
 with the scattered electric field $E^s(x,z) $, the wavenumber in vacuum $k^2_0 = \omega^2 \mu_0 \varepsilon_0$, the relative permittivity of the background medium $\varepsilon_{rb}$, and  the Green function $G(x,z|x',z')$ corresponding to the Helmholtz equation. The contrast function $\chi(x,z)$ is defined as:
 \begin{equation}
 \chi (x,z) = \frac{\varepsilon_r(x,z)}{\varepsilon_{rb}} - 1,
\end{equation}
 where $\varepsilon_r(x,z)$ denotes the relative permittivity containing the scattering object.
 
 For the sake of simplicity, we assume $\varepsilon_{rb} = 1$.
 
\subsection{Spatial discretization}
The same discretization in $x$ and $z$ as in \cite{Dilz2016The} and \cite{RJDilz2017spatialspectral} is employed here. Below, we give a brief summary of the discretization.


Along the $z$ direction, piecewise-linear functions $\Lambda_k$ are employed as expansion functions. The definition of $\Lambda_k$ is given by
\begin{small}
\begin{equation}
\Lambda_k (z) = \left \{
\begin{array}{ll} 
1 - \frac{|z - k\Delta - z_{min}|}{\Delta} \qquad \text{if} \qquad |z - k\Delta - z_{min}| < \Delta , \\
0 \qquad \text{if} \qquad |z  - k\Delta - z_{min}| > \Delta.
\end{array}\right. 
\end{equation}
\end{small}
Along the $x$ direction, a Gabor frame is employed to approximate $E(x,z)$ for fixed $z$. We have
\begin{small}
\begin{subequations}
\begin{equation}
\begin{split}
\chi (x, z_l) E^s (x, z_l)& = k^2_0  \chi(x,z_l)\int \limits^{z_{max}}_{z_{min}} \int \limits^{\infty}_{-\infty} G(x,z_l|x',z')\chi(x',z') E(x',z') dx' dz' \\
\end{split}
\end{equation}
\begin{equation}
\chi(x,z)E(x,z) = \chi E^i + \chi E^s = \sum^{N_k}_{k=0}\sum_{m,n} J_{mn,k}{g^w_{mn}}(x)\Lambda_k(z), 
\end{equation}
\label{chi_E}
\end{subequations}
\end{small}

where $N_k$ is the total number of expansion functions in the $z$ direction and 
\begin{small}
\begin{equation}
g^w_{mn}(x) = g^w(x - \alpha m X)e^{j\beta K n x}, \qquad {g^w}(x) = 2^{\frac{1}{4}}\exp \left( -\pi\frac{x^2}{X^2}\right),
\end{equation}
\end{small}
with the spatial window width $X$ and $K = 2\pi/X$. $\alpha$ and $\beta$ are the spatial and spectral oversampling rates respectively, such that $\alpha \beta  < 1$ . 
It is convenient to constrain $\alpha \beta$ to be a rational number in practice \cite{Bastiaans1995Gabor}, that is, $\alpha \beta =\frac{q}{p}$ with $q,p \in \mathbb{Z}^{+}$. And according to the calculation process of {the} Gabor dual window based on Zak transform, we need to choose smaller $q$ and $p$ to reduce the computation. In the following, we will find that the computation has no relation with $p$ and $q$ in our proposed method, thus we can choose oversampling parameters more freely.

According to the definition of the dual Gabor window, if $f(x)$ can be represented as
\begin{small}
$$ f(x) = \sum^{\infty}_{m =-\infty} \sum^{\infty}_{n =-\infty}f_{mn} {g^w_{mn}}(x),$$
\end{small}
then, the Gabor coefficients $\{f_{mn}\}$ can be obtained via
\begin{small}
\begin{equation}
f_{mn} = \int \limits^{\infty}_{-\infty} f(x) \eta^*_{mn}(x) dx ,
\label{gabor_coefficient_via_eta}
 \end{equation}
 \end{small}
where $ \eta_{mn}(x) = \eta (x - \alpha m X)e^{j\beta K n x}$ is the dual window and the asterisk denotes the conjugation operation.
\section{The Ewald Green function transformation}
\label{the_ewald_transformation}
For the Helmholtz equation in $\mathbb{R}^2$, i.e.
 \begin{small}
 \begin{equation}
 \Delta \phi + k^2_0 \phi = f,
 \end{equation}
 \end{small}
 the Green function as the fundamental solution that satisfies the radiation conditions at infinity for the time convention $\exp(j\omega t)$ of  the above equation is  
 \begin{small}
 \begin{equation}
G(x,z|x',z') = \frac{1}{4j}H^{(2)}_0 (k_0R),
 \end{equation}
 \end{small}
 where $H^{(2)}_0$ denotes the zeroth order Hankel function of the second kind \cite{Capolino2005Efficient} and  $R = \sqrt{(x-x')^2 + (z - z')^2}$.

\begin{figure}[!h]
  \centering\includegraphics[width=3.2in]{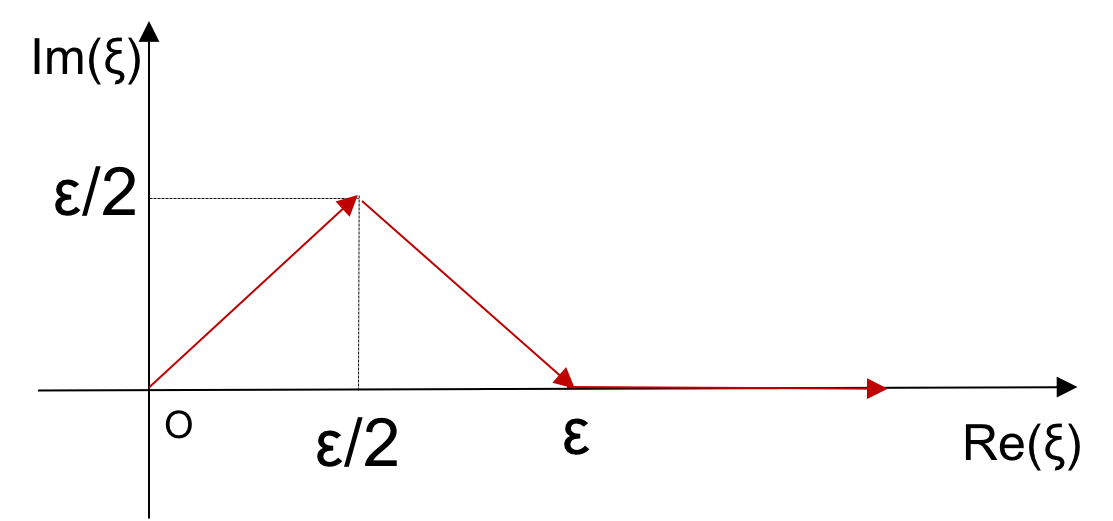}
 \caption{Contour path we employed in the Ewald transformation.}
 \label{integral_path}
 \end{figure}

 From \cite{Capolino2005Efficient}, we also obtain the Ewald transformation of the above Green function,
 \begin{small}
 \begin{equation}
 G(x,z|x',z')  = \frac{1}{2\pi}\int \limits^{\infty}_0\; \frac{\exp \left( -R^2 \xi^2 + \frac{k_0^2}{4\xi^2} \right)}{\xi} \; d\xi .
 \label{green_ewald}
 \end{equation}
 \end{small}
 According to the constraint conditions for the integration path of the Ewald transformation described in \cite{Capolino2005Efficient,Arens2013Analysing}, we utilize the integration contour, as shown in Fig.~\ref{integral_path}.

 The definition of the integration path is 
 \begin{small}
 \begin{equation}
 \xi(w) = \left\{\begin{array}{ll}
 (1 + j)w \qquad 0 <w< \mathcal{E}/{2}, \\
w + (\mathcal{E} - w) j \qquad \mathcal{E}/{2} < w< \mathcal{E}, \\
w  \qquad w > \mathcal{E} ,\end{array} \right. 
 \label{integral_path_p}
 \end{equation}   
 \end{small}
 where $ w \in (0, \infty)$. 

In the Ewald transformation, the integral in Eq.~\eqref{green_ewald} is split in two parts and we denote them as 
\begin{small}
\begin{subequations}
\begin{equation}
G_{spectral}(x,z|x',z') = \frac{1}{2\pi}\int \limits^{\mathcal{E}}_0 
\; \frac{\exp \left( -R^2\xi^2 + \frac{k^2_0}{4\xi^2} \right)}{\xi} \; d\xi, \\
\end{equation}
\begin{equation}
 G_{spatial}(x,z|x',z') = \frac{1}{2\pi}  \int \limits^{\infty}_\mathcal{E}
 \; \frac{\exp \left( -R^2\xi^2 + \frac{k^2_0}{4\xi^2} \right)}{\xi} \; d\xi .
 \end{equation}
 \end{subequations}
 \end{small}
\subsection{The spatial part of the Green function}
The spatial part of the Green function is
\begin{small}
\begin{equation}
G_{spatial}(x,z|x',z') = \frac{1}{2\pi} \int \limits^{\infty}_\mathcal{E} 
\;\frac{\exp \left( -R^2\xi^2 + \frac{k^2_0}{4\xi^2} \right)}{\xi} \; d\xi,
\end{equation}
\end{small}
where $\xi \in \mathbb{R}$ and $\xi \geq \mathcal{E}$. That means the integral path in the above equation is the real axis from $\mathcal{E}$ to positive infinity. The coordinates are now separable, owing to the square of $R$.
\subsection{The spectral part of the Green function}
The spectral part of the Green function is
\begin{small}
\begin{equation}
G_{spectral}(x,z|x',z') = \frac{1}{2\pi} \int \limits^{\mathcal{E}}_0 
\;\frac{\exp \left( -R^2\xi^2 + \frac{k^2_0}{4\xi^2} \right)}{\xi} \; d\xi .
\label{G_spectral_2.2}
\end{equation}
\end{small}
Its Fourier transform is 
\begin{small}
\begin{equation}
\label{G_spectral_original}
\widehat{G}_{spectral}(k_x,z|z')
= \frac{1}{2\pi} \int \limits^{\mathcal{E}}_0  \;\frac{1}{\xi} e^{-(z-z')^2\xi^2 + \frac{k_0^2}{4\xi^2}} \left[  \; \int \limits^{\infty}_{-\infty} \; e^{-x^2 \xi^2 - jk_x x} \; dx\right] d\xi.
\end{equation}
\end{small}
We require that $\mathrm{Re}\{\xi^2\} > 0$ to guarantee the convergence of the $\xi$-integral in Eq.~\eqref{G_spectral_2.2}, which is automatically satisfied for the chosen path as shown in Fig.~\ref{integral_path}. Once this is ensured, the result is not affected by the interchange of the order of integration. According to  the closed-form evaluation of the $\xi$-integral using the formula (2.33) in \cite{2014Table}
\begin{small}
 \begin{equation}
\int \limits^{\infty}_{-\infty} e^{-a\xi^2 + b\xi}\; d\xi = \sqrt{\frac{\pi}{a}}e^{\frac{b^2}{4a}},
 \label{xi_integral}
 \end{equation}
 \end{small}
 we get
 \begin{small}
\begin{equation}
\widehat{G}_{spectral}(k_x,z|z') \xlongequal{\zeta = 1/\xi} \frac{1}{2\sqrt{\pi}}\int \limits^{\infty}_{1/\mathcal{E}} \;  e^{-(z-z')^2/\zeta^2+ (k_0^2 - k_x^2)\zeta^2/4} \; d\zeta ,
\end{equation}
\end{small}
where the integral path $\xi$ is defined by Eq.~\eqref{integral_path_p} and
$ \zeta = \frac{1}{\xi(w)}$. Replacing $w$ by $\frac{1}{w}$ does not change the path $\zeta$, thus we rewrite $\zeta$ as
\begin{small}
\begin{equation}
\zeta 
=\frac{1}{\xi(\frac{1}{w}) } 
=
 \left\{ \begin{array}{ll}
\frac{w - (\mathcal{E}w^2 -w)j}{1 + (\mathcal{E}w - 1)^2} \qquad {1}/\mathcal{E} < w < {2}/\mathcal{E}, \\
\\
\frac{1-j}{2} w \qquad {2}/{\mathcal{E}} < w < \infty .\\
\end{array}\right.  
\label{integral_path_q}
\end{equation}
\end{small}
\section{Calculation of the domain integral representation}
 Substituting the Ewald transformation Eq.~\eqref{green_ewald} in Eq.~\eqref{chi_E}, we obtain
 \begin{small}
 \begin{equation}
\begin{split}
 &\chi(x,z_l) E^s(x,z_l)  
 = \frac{k^2_0}{2\pi} \sum_{m,n,k}  J_{mn,k}\cdot  \\
& \quad \int \limits^{\infty}_0    \; 
 \frac{ e^{\frac{k_0^2}{4\xi^2}}}{\xi}  \cdot 
  \left[ \; \int \limits^{\infty}_{-\infty} dx' \; e^{-(x-x')^2\xi^2} \cdot  {g^w_{mn}}(x') \right]\cdot  
 \left[\; \int \limits^{z_{max}}_{z_{min}} dz' \; \Lambda_k(z')\cdot e^{-(z' - z_l)^2\xi^2} \right] d\xi
\end{split}
 \end{equation}
 \end{small}
Owing to the Ewald transformation, the integrals related to $x'$ and $z'$ are separable. We split the above integral in two parts,
 \begin{small}
 \begin{subequations}
 \begin{equation}
\begin{split}
\chi E_{spatial}(x,z_k) &= \int \limits^{z_{max}}_{z_{min}} \int \limits^{\infty}_{-\infty} G_{spatial}(x,z_k|x',z') \cdot \chi E(x',z')\; dx'dz' \\
&= \sum_{m,n} f^{spatial}_{mn,k} {g^w_{mn}}(x),
\end{split}
\label{chi_E_spatial}
\end{equation}
\begin{equation}
\begin{split}
\chi E_{spectral}(x,z_k) &= \int \limits^{z_{max}}_{z_{min}} \int \limits^{\infty}_{-\infty} G_{spectral}(x,z_k|x',z') \cdot \chi E(x',z')\; dx'dz' \\
&=  \sum_{m,n} f^{spectral}_{mn,k} {g^w_{mn}}(x).
\end{split}
\label{chi_E_spectral}
 \end{equation}
 \end{subequations}
 \end{small}
In the following, we further simplify these two representations by working out the integrals over $x'$ and $z'$.
\subsection{Spatial part of the electric field}
According to Eq.~\eqref{chi_E_spatial} and Eq.~\eqref{gabor_coefficient_via_eta}, we have
\begin{small}
\begin{equation}
\begin{split}
f^{spatial}_{st,l} &= \int \limits^{\infty}_{-\infty} \chi E_{spatial}(x,z_l) \cdot \eta^*_{st}(x) dx \\
&= \frac{1}{2\pi}\sum_{m,n,k}J_{mn,k} \int \limits^{\infty}_{\mathcal{E}}  \; \frac{e^{\frac{k_0^2}{4\xi^2}}}{\xi} \cdot \left[ \; \int \limits^{z_{max}}_{z_{min}}\; \Lambda_k(z') \cdot e^{-(z_l-z')^2\xi^2 } \;  dz'\right]
\cdot \\
&\qquad  \quad 
 \Bigg{\{}  \int \limits^{\infty}_{-\infty}  \; \eta^*_{st}(x) \left[\; \int \limits^{\infty}_{-\infty} \; e^{-(x-x')^2\xi^2} \cdot {g^w_{mn}}(x') \; dx'\right] \; dx\Bigg{\}} d\xi .
\end{split}
\label{f_spatial_mnl}
\end{equation}
\end{small}

\subsubsection{The integral related to $x$ and $x'$}
\label{integral_x_spatial}
The integral related to  $x$ and $x'$ in Eq.~\eqref{f_spatial_mnl} is 
\begin{small}
\begin{equation}
\begin{split}
& \quad  \int \limits^{\infty}_{-\infty} \; \eta^*_{st}(x) \left[ \; \int \limits^{\infty}_{-\infty} \; e^{-(x-x')^2\xi^2} \cdot g_{mn}(x') \; dx'\right] \; dx \\
&=  \int \limits^{\infty}_{-\infty} \int \limits^{\infty}_{-\infty}  \; e^{-(x-x')^2\xi^2} \cdot {g^w}(x' - \alpha mX) \cdot e^{j\beta K n x'} \cdot \eta^* (x - \alpha s X) \cdot e^{-j\beta K t x} \; dx'dx .\\
\end{split}
\label{integral_eta_g_e}
\end{equation}
\end{small}
We assume that the dual Gabor window $\eta(x)$ can be represented by a weighted sum of modulated Gaussian functions of the Gabor frame used for the expansion as follows,
\begin{equation}
\eta (x) = \sum^{N_u}_{u} \sum^{N_v}_{v} a_{u v}\cdot {g^w_{uv}}(x).
\end{equation}
By substituting the above equation in Eq.~\eqref{integral_eta_g_e}, we obtain
\begin{small}
\begin{equation}
\begin{split}
& \quad  \int \limits^{\infty}_{-\infty} \; \eta^*_{st}(x) \left[ \;  \int \limits^{\infty}_{-\infty} \; e^{-(x-x')^2\xi^2} \cdot {g^w_{mn}}(x') \; dx'\right]dx \\
&=\sqrt{2}\sum^{N_{u}}_{u} \sum^{N_{v}}_{v} a^*_{u v} \cdot e^{j2\pi\alpha \beta s v} \cdot  \\
&\qquad \int \limits^{\infty}_{-\infty} \int \limits^{\infty}_{-\infty} \; e^{-(x-x')^2\xi^2 - \frac{\pi}{X^2} (x' - \alpha m X)^2 - \frac{\pi}{X^2}(x - \alpha s X - \alpha u X)^2 + j\beta K n x' -j\beta K (v + t) x } \; dx'  dx\\
&\xlongequal[q = m-s-u]{p = v + n + t}\sum^{N_{u}}_{u} \sum^{N_{v}}_{v}2\sqrt{\pi}X^2 a^*_{u v} \cdot e^{j2\pi\alpha \beta [s v + mn - (s+u)(t+v)] - \frac{\pi}{2}\beta^2(v + t -n)^2} \cdot f(q,p,\xi),
\end{split}
\end{equation}
\end{small}
where
\begin{small}
\begin{equation}
 f (q,p,\xi) = \frac{1}{\sqrt{4X^2\xi^2 + 2\pi}} \cdot \exp \left( - \frac{\pi}{2 + \frac{\pi}{X^2\xi^2}} \cdot (\alpha q + j\beta p)^2 - \frac{\pi}{2} \beta^2p^2  \right).
\end{equation}
 \end{small}
The term $- \frac{\pi}{2}\beta^2 p^2 $ in the argument of the exponential function ensures that the real part of the  exponent 
\begin{small}
 $$ - \frac{\pi}{2 + \frac{\pi}{X^2\xi^2}} \cdot (\alpha q + j\beta p)^2 - \frac{\pi}{2}\beta^2 p^2 $$
 \end{small}
 is always negative. Thus, for any fixed $q$ and $p$, $f(q,p,\xi)$ 
remains bounded.

 \subsubsection{The integral related to $z'$}
 \label{integral_z_spatial}
 After mathematical deduction we obtain the analytical results of the integral related to $z'$. Let
 \begin{small}
 \begin{equation}
 \begin{split}
 g({k-l},\xi) &= ({k-l}+1)\cdot \frac{\sqrt{\pi}}{2\xi} \left[ \text{erf}(\xi ({k-l}+1)\Delta) - \text{erf}(\xi {k-l} \Delta) \right]  +  \\
& \quad  \frac{1}{2\Delta \xi^2} \left[ e^{-({k-l}+1)^2\Delta^2 \xi^2} - e^{-({k-l})^2\Delta^2 \xi^2} \right], 
\end{split}
 \end{equation}
 \end{small}
 where $\text{erf} (z)$ is the error function defined as \cite{M1972Handbook}
 \begin{equation}
 \text{erf}(\xi) = \frac{2}{\sqrt{\pi}} \int \limits^\xi_{-\infty} e^{-t^2} dt.
 \end{equation}
 Then we denote
 \begin{small}
 \begin{equation}
 h(k-l,\xi)  \triangleq \int \limits^{z_{max}}_{z_{min}} \Lambda_k(z') \cdot e^{-(z_l-z')^2\xi^2}dz' =
 \left\{ \begin{array}{ll}
 g(k-l,\xi) \quad k = 0, \\ 
 g(k-l,\xi) + g(l-k,\xi) \quad 0 < k < N_k, \\ 
 g(l-k,\xi) \quad k = N_k .
 \end{array}\right.
 \end{equation}
 \end{small}

 \subsubsection{The entire integral of the spatial part}
We denote
 \begin{small}
 \begin{equation}
 \begin{split}
 &\quad f^{spatial}_{st,l}(m,n,k) \\
 &=  \frac{1}{2\pi} \int \limits^{\infty}_{\mathcal{E}}  \frac{e^{\frac{k^2_0}{4\xi^2}}}{\xi} \cdot \Bigg{\{} \; \int \limits^{\infty}_{-\infty}  \; \eta^*_{st}(x') \left[ \; \int \limits^{\infty}_{-\infty} \; e^{-(x-x')^2\xi^2} \cdot g_{mn}(x) dx \right] dx' \Bigg{\}} \cdot  \\
 &\qquad  \quad \left[ \; \int \limits^{z_{max}}_{z_{min}} \; \Lambda_k(z') \cdot e^{-(z_l-z')^2\xi^2}  dz'\right] d\xi \\ 
&= \frac{X^2}{\sqrt{\pi}} \sum_{u,v} a^*_{u v} \cdot e^{j2\pi\alpha\beta[sv +mn - (s+u)(t+v)] - \frac{\pi}{2} \beta^2 (v + n -n)^2}\cdot  \\
 &\qquad \quad \int \limits^{\infty}_{\mathcal{E}} \frac{e^{\frac{k^2_0}{4\xi^2}}}{\xi} \cdot f(m-s-u, n + t + v,\xi) \cdot h(k-l,\xi)\; d\xi . \\
 \end{split}
 \label{f_spatial_mnl_stk}
 \end{equation}
 \end{small}
  The definitions of $f(q,p,\xi)$ and $h(k-l,\xi)$ can be found in Section \ref{integral_x_spatial} and \ref{integral_z_spatial}, respectively. Consider the following integral :
  \begin{small}
  \begin{equation}
 \begin{split}
 & \quad \int \limits^{\infty}_{\mathcal{E}} d\xi  \cdot \frac{e^{\frac{k^2_0}{4\xi^2}}}{\xi} \cdot f(q, p,\xi) \cdot h(k-l ,\xi)\; d\xi  \\
 &= h_1 \int \limits^{\infty}_{\mathcal{E}}\frac{e^{\frac{k^2_0}{4\xi^2}}}{\xi} \cdot f(q, p,\xi) \cdot  g(k-l,\xi) \;  d\xi  +
 h_2 \int \limits^{\infty}_{\mathcal{E}} \frac{e^{\frac{k^2_0}{4\xi^2}}}{\xi} \cdot f(q, p,\xi) \cdot  g(l-k ,\xi) \; d\xi,\\
 \end{split}
 \end{equation}
 \end{small}
 where $h_1, h_2 \in \{0,1\}$. Thus, we can study the following alternative \mbox{integral:}
 \begin{small}
 \begin{equation}
   \int \limits^{\infty}_{\mathcal{E}} \;\frac{e^{\frac{k^2_0}{4\xi^2}}}{\xi} \cdot f(q, p,\xi) \cdot  g(k-l ,\xi) \; d\xi .
   \label{integral_spatial}
 \end{equation}
 \end{small}
 If we want to calculate $f^{spatial}_{st,l}(m,n,k)$ for every $m,s \in \{-N_m,\ldots ,N_m\}$, $n,t \in \{-N_n,\ldots,N_n\}$ and  $k,l \in \{0,\ldots, N_k\}$, we need to calculate Eq.~\eqref{integral_spatial}
 for every $q$, $p$ and $k-l$. 

\subsection{Spectral part of electric field}
 Based on the convolution property, the Fourier transform of the spectral part, i.e. Eq.~\eqref{chi_E_spectral} can be simplified as
 \begin{small}
 \begin{equation}
 \begin{split}
 \widehat{\chi E}_{spectral}(k_{x},z_l) &=  \sum_{m,n,k}J_{mn,k} \mathcal{F} \left(  \; \int \limits^{z_{max}}_{z_{min}}\int \limits^{\infty}_{-\infty} G_{spectral}(x,z_l|x',z')\cdot {g^w_{mn}}(x') \cdot \Lambda_k (z')  dx' dz'  \right) \\
 &=\frac{1}{2\sqrt{\pi}} \sum_{m,n,k}J_{mn,k}\cdot {\widehat{g^w}_{nm}}(k_{x}) \cdot e^{2\pi j \alpha \beta mn} \cdot \\
 &\qquad \quad \int \limits^{\infty}_{1/\mathcal{E}}\;  e^{(k_0^2 - k_{x}^2)\zeta^2/4} \cdot  \left[ \; \int \limits^{z_{max}}_{z_{min}} \; \Lambda_k(z') \cdot e^{-(z_l - z')^2/\zeta^2} dz'\right] d\zeta,
 \end{split}
 \end{equation}
 \end{small}
 where $\zeta = q(w)$ and the definition of $q(w)$ is shown in \eqref{integral_path_q}.
Let
 \begin{small}
 \begin{equation}
  \widehat{\chi E}_{spectral}(k_{x},z') = \sum^{N_k}_{k=0} \sum_{n,m} \widetilde{E}_{nm,k}{\widehat{g^w}_{nm}}(k_{x})\Lambda_k(z'),
 \end{equation}
 \end{small}
 then we have
 \begin{small}
 \begin{equation}
 \begin{split}
 \widetilde{E}_{ts,l} &= \int^{\infty}_{-\infty} \widehat{\chi E}_{spectral}(k_{x},z_l) \cdot \widehat{\eta}^*_{ts}(k_{x}) dk_{x}\\
 &= \frac{1}{2\sqrt{\pi}}\sum^{N_k}_{k=0}\sum_{m,n}J_{mn,k}\cdot e^{2\pi j \alpha \beta mn} 
 \int \limits^{\infty}_{1/\mathcal{E}} \; e^{k_0^2 \zeta^2/4}\cdot \left[ \; \int \limits^{z_{max}}_{z_{min}}\; \Lambda_k(z') \cdot e^{-(z_l- z')^2/\zeta^2}dz' \right]
  \cdot  \\
&\qquad \quad \left[\;  \int \limits^{\infty}_{-\infty} e^{-k^2_{x}\zeta^2/4}\cdot {\widehat{g^w}_{nm}}(k_{x})\cdot \widehat{\eta}^*_{ts}(k_{x})dk_{x}\right]d\zeta 
\end{split}
\label{E_spectral_nml}
\end{equation}
\end{small}
where 
\begin{small}
\begin{equation}
\begin{array}{ll}
{\widehat{g^w}} (k_{x}) = 2^{\frac{1}{4}} X e^{- \frac{\pi}{K^2}k^2_{x}}, \qquad  
{\widehat{g^w}_{nm}}(k_{x}) ={ \widehat{g^w}}(k_{x} - n\beta K)e^{-jm\alpha Xk_{x}}, \\
\widehat{\eta}(k_{x}) = \sum \limits^{N_{u}}_{u}  \sum \limits^{N_{v}}_{v} \widehat{a}_{u v} {\widehat{g^w}_{u v}}(k_{x}),
 \qquad
\widehat{\eta}_{ts}(k_{x}) = \widehat{\eta}(k_{x} - t\beta K)e^{-js\alpha X k_{x}}.
\end{array}
\end{equation}
\end{small}
\subsubsection{The integral related to $k_{x}$}
\label{integral_kx_spectral}
 The integral related to $k_{x}$ in Eq.~\eqref{E_spectral_nml} can be solved analytically as follows, 
 \begin{small}
 \begin{equation}
 \begin{split}
 &\quad   \int \limits^{\infty}_{-\infty} dk_{x} \; e^{-k^2_{x}\zeta^2/4} \cdot  {\widehat{g^w}_{nm}}(k_{x}) \cdot \widehat{\eta}^*_{ts}(k_{x})  \\
 &\triangleq 2^{\frac{3}{2}}X^2K \sum^{N_{u}}_{u} \sum^{N_{v}}_{v} \widehat{a}^*_{uv} \cdot e^{-j2\pi \alpha \beta t v} \cdot e^{-\frac{\pi}{2}\beta^2 (n-t -u)^2} \cdot \widetilde{f}(q,p,\zeta),
 \end{split}
 \end{equation}
 \end{small}
 where
 \begin{small}
 \begin{equation}
 \left\{
 \begin{array}{ll}
 q = s + v - m ,\\
 p = t + u + n , \\
 \widetilde{f}(q, p, \zeta) = \sqrt{ \frac{\pi}{K^2\zeta^2 + 8\pi}}\cdot  \exp \left( \frac{4\pi^2}{K^2\zeta^2+8\pi}(\beta p + j\alpha q)^2 - \frac{\pi}{2}\beta^2p^2 \right) .
 \end{array}\right.
 \end{equation}
 \end{small}

\subsubsection{The integral related to $z'$}
\label{integral_z_spectral}
 Similar to the spatial part, we have
 \begin{small}
 \begin{equation}
 \widetilde{h}(k-l,\zeta) \triangleq \int \limits^{z_{max}}_{z_{min}}\; \Lambda_k(z') \cdot e^{-( z_l - z')^2/\zeta^2} \; dz'
 = \left\{ \begin{array}{ll} 
\widetilde{g}(k-l,\zeta) \quad k = 0,\\
\widetilde{g}(k-l,\zeta) + \widetilde{g}(l-k,\zeta) \quad 0 < k< N_k ,\\
\widetilde{g}(l-k,\zeta) \quad k = N_k ,
\end{array}\right.
\end{equation}
\end{small}
and 
\begin{small}
\begin{equation}
\begin{split}
\widetilde{g}(k-l,\zeta) &= \frac{\sqrt{\pi}(k-l+1) }{2} \zeta \cdot \left[ \text{erf}\left( \frac{(k-l+1)\Delta}{\zeta} \right) - \text{erf}\left( \frac{(k-l)\Delta}{\zeta} \right) \right] + \\&\qquad \quad \frac{\zeta^2}{2\Delta}\cdot  \left[ \exp \left( -\frac{(k-l+1)^2\Delta^2}{\zeta^2} \right) - \exp \left( -\frac{(k-l)^2\Delta^2}{\zeta^2} \right) \right] .
\end{split}
\label{widetilde_g}
\end{equation}
\end{small}

\subsubsection{The entire integral of the spectral part}
 We denote
 \begin{small}
 \begin{equation}
 \begin{split}
 \widetilde{E}_{ts,l}(m,n,k) 
 &= \frac{1}{2\sqrt{\pi}}\cdot e^{2\pi j \alpha \beta mn} 
\int \limits^{\infty}_{1/\mathcal{E}} \; e^{k_0^2 \zeta^2/4} \left[ \; \int \limits^{z_{max}}_{z_{min}}\; \Lambda_k(z) \cdot e^{-(z_l- z')^2/\zeta^2} dz'\right]\ \cdot  \\
 &\qquad  \qquad \left[ \int \limits^{\infty}_{-\infty} \; e^{-k^2_{x}\zeta^2/4} \cdot  {\widehat{g^w}_{nm}}(k_{x}) \cdot \widehat{\eta}^*_{ts}(k_{x})dk_{x} \right] d\zeta \\
 &= 2X\sqrt{2\pi}\sum_{uv}\widehat{a}^*_{uv}\cdot e^{j2\pi\alpha\beta(mn - tv) -\frac{\pi}{2}\beta^2(n-t-u)^2}\cdot \\
 &\quad \int \limits^{\infty}_{1/\mathcal{E}} e^{k_0^2\zeta^2/4} \cdot \widetilde{f}(s+v-m,t+u+n,\zeta) \cdot \widetilde{h}(k-l,\zeta) d\zeta .
 \end{split}
 \label{E_spectral_nml_stk}
 \end{equation}
 \end{small}
 The definitions of $\widetilde{f}(q,p,\zeta)$ and $\widetilde{h}(k-l, \zeta)$ can be found in Section \ref{integral_kx_spectral} and Section \ref{integral_z_spectral}, respectively.

 Considering that $\widetilde{h}(k-l,\zeta)$ can be directly expressed in $\widetilde{g}(k-l,\zeta)$, we can study the following alternative integral :
 \begin{small}
 \begin{equation}
\int \limits^{\infty}_{1/\mathcal{E}} \; e^{k_0^2\zeta^2/4} \cdot \widetilde{f}(q,p,\zeta) \cdot \widetilde{g}(k-l,\zeta) \; d\zeta.
 \label{integral_spectral}
 \end{equation}
 \end{small}
\section{Splitting parameter $\mathcal{E}$}
\label{optimum_splitting_parameter}
\subsection{Asymptotic {c}onvergence of the integrands}
Here we analyze the convergence properties of the integrals related to $G_{spatial}$ and $G_{spectral}$ given by Eq.~\eqref{integral_spatial} and Eq.~\eqref{integral_spectral}, respectively.

\subsubsection{Spatial part}
 For large $\xi$, it is easy to obtain the asymptotic convergences of $\frac{e^{\frac{k^2_0}{4\xi^2}}}{\xi}$ {and} $f(q,p,\xi)${,} respectively,
\begin{small}
 \begin{equation}
\frac{e^{\frac{k^2_0}{4\xi^2}}}{\xi} \sim \frac{1}{\xi},
 \end{equation}
\begin{equation}
 \big| f(q,p,\xi) \big| \sim  \frac{1}{\xi} \cdot e^{ - \frac{\pi}{2}\alpha^2 p^2}.
 \end{equation}
\end{small}
From \cite{Capolino2005Efficient}, we obtain the asymptotic expansion for the error function for large argument, $\text{erf}(\xi) \sim 1 - \frac{e^{-\xi^2}}{\sqrt{\pi}\xi}$. Therefore the asymptotic behavior of $g(k-l,\xi)$ for large $\xi$ {is}:
 \begin{small}
 \begin{equation}
 g(k-l,\xi) \sim \left\{ \begin{array}{ll}
 \frac{1}{2\Delta \xi^2} \qquad k = l-1, \\
 \frac{\sqrt{\pi}}{2\xi} \qquad   k = l, \\
 \frac{1}{2(k-l) \Delta \xi^2}e^{-(k-l)^2\Delta^2 \xi^2} \qquad \text{other}.
 \end{array}\right.
 \end{equation}
 \end{small} 
 As a consequence, we have
 \begin{small}
 \begin{equation}
 \frac{e^{\frac{k_0^2}{4\xi^2}}}{\xi} \cdot f(q,p,\xi) \cdot g(k-l,\xi) \sim 
 \left\{\begin{array}{ll}
  \frac{1}{4X\Delta \xi^4} \cdot \exp \left( -\frac{\pi \alpha^2 q^2}{2} \right) \quad  k=l -1 ,\\
 \frac{\sqrt{\pi}}{4X\xi^3} \cdot \exp \left( -\frac{\pi \alpha^2 q^2}{2} \right) \quad  k= l ,\\
 \frac{1}{4(k-l)X\Delta \xi^4} \cdot \exp \left( - (k-l)^2\Delta^2 \xi^2 - \frac{\pi \alpha^2 q^2}{2} \right) \quad \text{other}.
 \end{array}\right.
 \label{limit_spatial}
 \end{equation}
 \end{small}
 when $\xi$ approaches infinity.
 If $\frac{\pi\alpha^2q^2}{2}$ is large enough, $\frac{e^{\frac{k_0^2}{4\xi^2}}}{\xi} \cdot f(q,p,\xi) \cdot g(k-l,\xi)$ is close to zero.

 \subsubsection{Spectral part}
 When $w$ approaches infinity $(> \frac{2}{\mathcal{E}})$, we have
 \begin{small}
 \begin{equation}
  \big| e^{-k^2_x\zeta(w)^2/4} \big|  = \big | \exp \left( -j k^2_x w^2/8 \right) \big| \sim 1
 \end{equation}
 \begin{equation}
 \big| \widetilde{f}(q, p, \zeta(w))  \big| \sim 
 \sqrt{2\pi}\cdot \frac{1}{Kw} \cdot \exp\left( -\frac{\pi}{2}\beta^2 p^2 \right).
 \end{equation}
\end{small}
By expanding the integrand $e^{-\zeta^2}$ into its Maclaurin series and integrating term by term, we obtain the error function's Maclaurin series as
\begin{small}
\begin{equation}
\text{erf}(\zeta) = \frac{2}{\sqrt{\pi}} \sum^{\infty}_{n=0} \frac{(-1)^n \zeta^{2n+1}}{n!(2n+1)}.
\label{err_Maclaurin}
\end{equation}
\end{small}
By substituting Eq.~\eqref{err_Maclaurin} {in} Eq.~\eqref{widetilde_g} and letting $w\rightarrow \infty$ (thus $\zeta$ {goes  to $(1-j)\infty$}), we can obtain the asymptotic behavior of $\widetilde{g}$ as follows,
\begin{small}
\begin{equation}
\widetilde{g}(k-l,\zeta) \sim \frac{\Delta}{2}.
\end{equation}
\end{small}

 As a consequence, when $w \rightarrow \infty$,
 \begin{small}
 \begin{equation}
 \big|  e^{k_0^2\zeta^2(w)/4} \cdot \widetilde{f}(q,p,\zeta(w)) \cdot \widetilde{g}(k-l,\zeta(w)) \cdot \zeta'(w) \big|
 \sim
 \frac{\sqrt{2\pi}\Delta (1 -j)}{4Kw} \cdot \exp \left( -\frac{\pi}{2}\beta^2 p^2\right) .
 \label{limit_spectral}
 \end{equation}
 \end{small}

\subsection{Discussion on the choice of $\mathcal{E}$}
 In the $spatial$ part, for large $\xi$,  the asymptotic behavior of the integrand exhibits Gaussian convergence $e^{-(k-l)^2\Delta^2\xi^2}$ if $k \neq l-1, l$. In the $spectral$ part, the large argument ${w}> \frac{2}{\mathcal{E}}$ leads to {a} rapidly oscillatory behavior in the integrand, $e^{k_0^2 \zeta^2(w)/4} = \exp \left( - j\frac{k^2_0 w^2}{8}\right)$. If $k^2_0 w^2$ is large enough, the infinite integration interval can be truncated. To minimize the calculation {while} holding high accuracy, we would like to maximize the lower bounds of $(k-l)^2\Delta^2\xi^2$ and $\frac{k^2_0 w^2}{8}$. Because
 \begin{small}
 \begin{equation}
 (k-l)^2\Delta^2\xi^2 \geq \Delta^2 \mathcal{E}^2  \qquad  \text{for} \; k \neq l-1, l\; \text{and}\;\xi \geq \mathcal{E}  ,
 \end{equation}
 \begin{equation}
 \frac{k^2_0 w^2}{8} \geq \frac{k^2_0}{2\mathcal{E}^2} \qquad \text{for} \; w \geq \frac{2}{\mathcal{E}}.
 \end{equation}
 \end{small}
 Thus, we let
 \begin{small}
 \begin{equation}
 \Delta^2 \mathcal{E}^2 = \frac{k^2_0}{2\mathcal{E}^2},
 \end{equation}
 \end{small}
 to get the optimum splitting parameter 
 \begin{small}
 \begin{equation}
 \mathcal{E} = 2^{-\frac{1}{4}}\sqrt{\frac{k_0}{\Delta}}.
 \end{equation}
 \end{small}
\section{Numerical examples}
\label{numerical_examples}
\subsection{Validation and accuracy}
We employ the same scatters described in \cite{Dilz2016The} to test the accuracy of our algorithm. The parameters of the three scatters are listed in Table \ref{scatters_parameters}.

$E^i(x,z) = E_0 \exp (jk_0 (x \cos \theta + z \sin \theta))$ with $E_0 = 1V/m$. {After several numerical experiments, we construct the approximation of the dual window by choosing $u=2$ and $v=3$.}

Furthermore, the parameters of the discretization in the $x$ and $z$ direction are shown in Table \ref{discretization_parameters}. As a consequence, there are $13$ ($|m| \leq M, m \in \mathbb{Z}$) spatial windows with $0.5m$ width and $7$ ($|n| \leq N, n \in \mathbb{Z}$) spectral windows.
\setcounter{table}{0}

\begin{small}
\begin{table}[h]
\centering
\begin{tabular}{|c|c|c|c|}
\hline
Parameters & Circle & Rectangle & Grating  \\\hline 
$\varepsilon_r$ of object & 2 & 2 & 2\\\hline
size of object (m) & $r = 1.35$ & $2.0 \times 5.0$ & 5 blocks, $1\times 1.4$, spacing $2$ \\\hline
wavenumber $k_0$ (m$^{-1}$) & 1.45 & 0.8388 & 1.5 \\\hline
angle of incidence $\theta$ & $0^{\circ}$ & $90^{\circ}$ & $45^{\circ}$ \\\hline
\end{tabular}
\caption{The parameters of the three {scattering} examples.}
\label{scatters_parameters}
\end{table}
\end{small}

\begin{small}
\begin{table}[h]
\centering
\begin{tabular}{|c|c|c|c|c|c|c|}
\hline
 Parameters of the discretization & $X$ & $M$ & $N$ & $\alpha (=\beta)$ & $\Delta $ & $N_k$ \\\hline
 $x$-direction &  $0.5$ & $6$ & $3$ & $\sqrt{2/3}$ &   \diagbox[dir=SW]{}{\textcolor{white}{.}}  &   \diagbox[dir=SW]{}{\textcolor{white}{.}} 
  \\\hline
  $z$-direction  &  \diagbox[dir=SW]{\textcolor{white}{.}}{} & \diagbox[dir=SW]{}{\textcolor{white}{.}}  &   \diagbox[dir=SW]{}{\textcolor{white}{.}}  &   \diagbox[dir=SW]{}{\textcolor{white}{.}}  & 0.05 & 56
\\\hline
\end{tabular}
\caption{The parameters of the discretization in the $x$ and $z$ direction.}
\label{discretization_parameters}
\end{table}
\end{small}

\begin{figure}[htb]
\centering 
\subfigure[]{ 
\label{example1_realpart} 
\includegraphics[height=4.5cm]{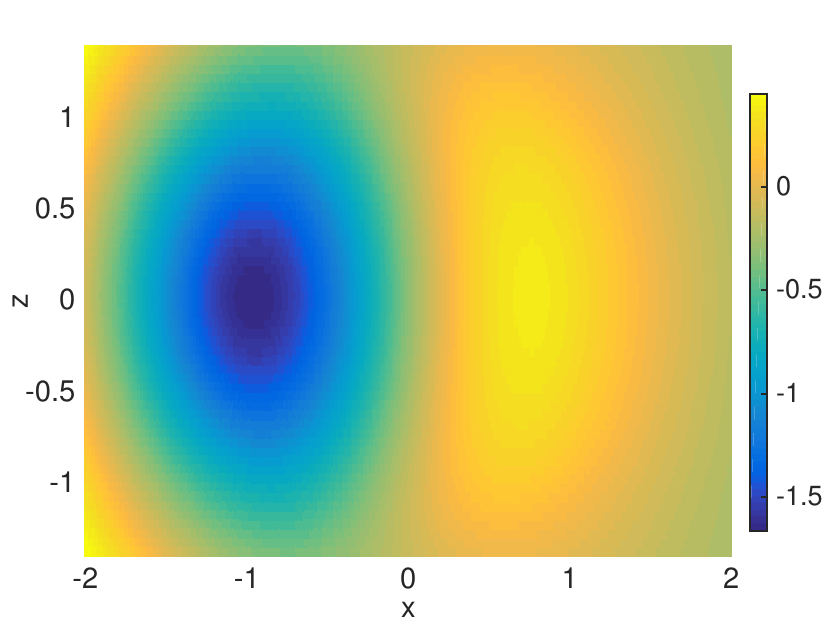}} 
\subfigure[]{ 
\label{example2_realpart}
\includegraphics[height=4.5cm]{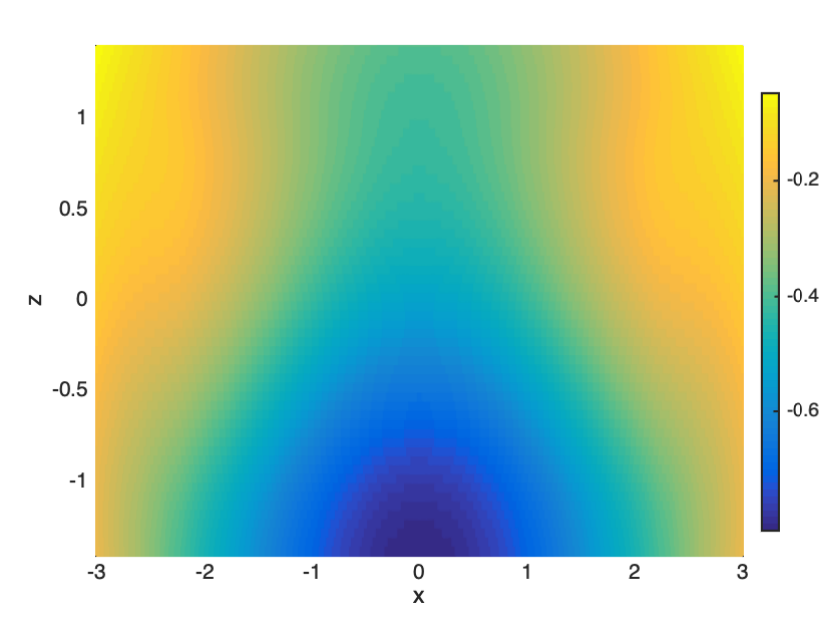}} 
\caption{The real parts of $\chi E^s (x,z)$ for the circular and rectangular scatterers.} 
\label{realpart}
\end{figure}

 \begin{figure}[htb]
\centering \subfigure[]{ 
\label{example1_error} 
\includegraphics[height=4cm]{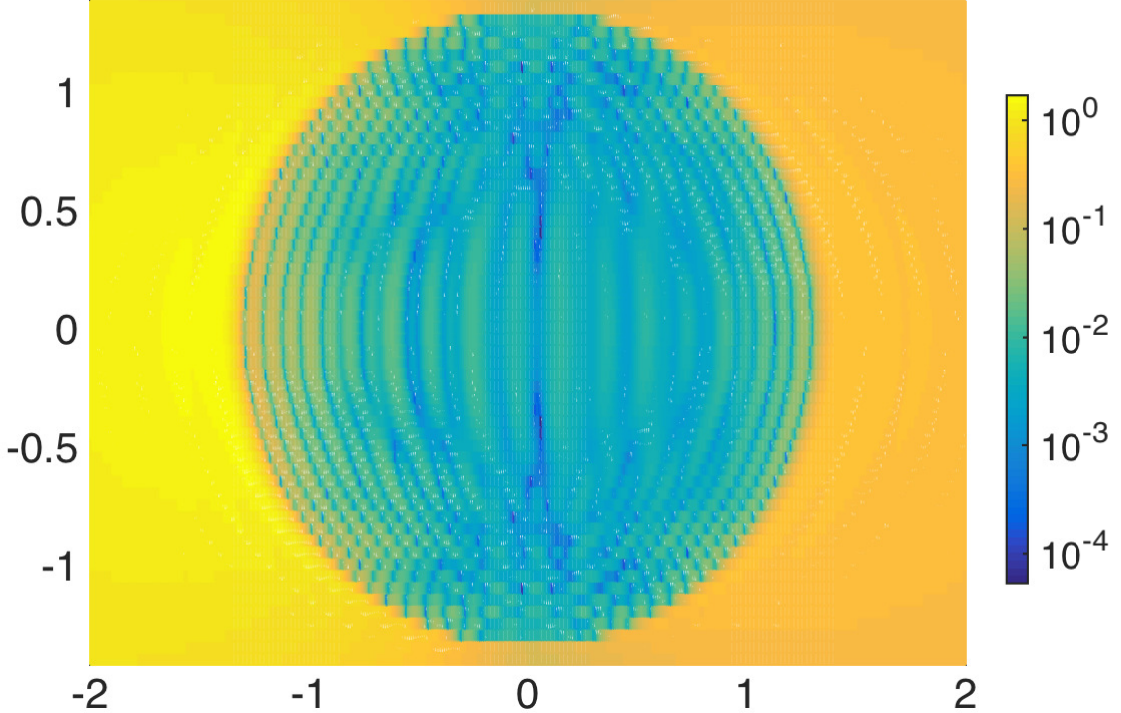}} 
\subfigure[]{ 
\label{example2_error}
 \includegraphics[height=4cm]{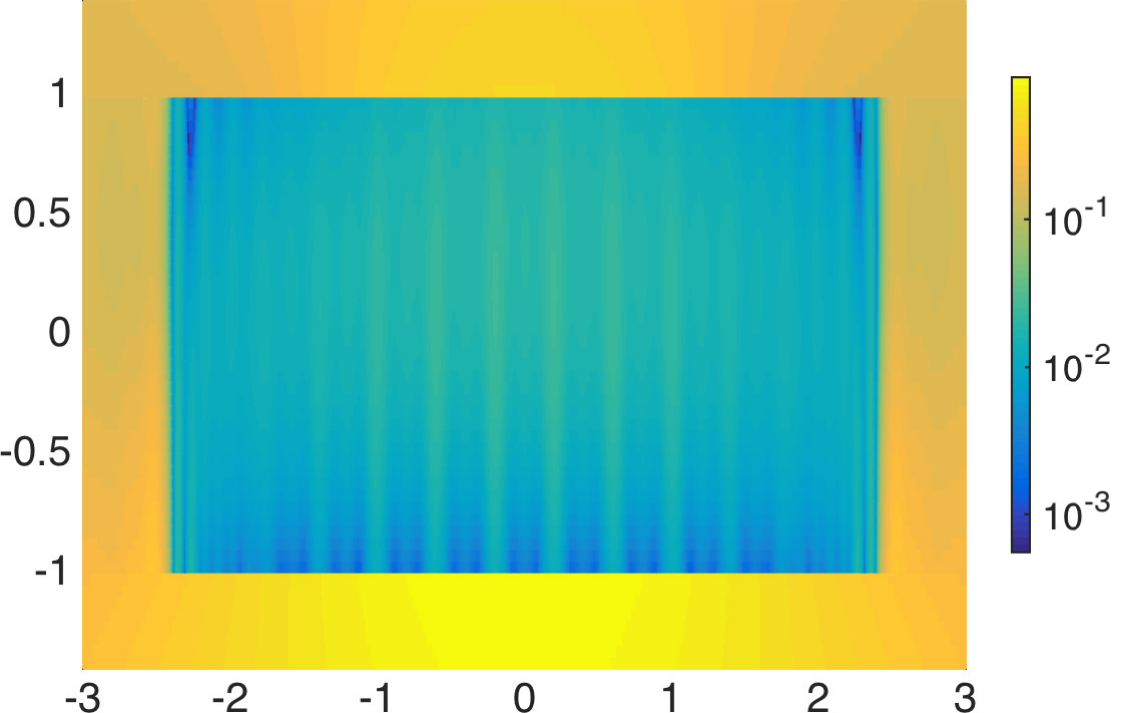}} 
\caption{The absolute value of the error between $\chi E^s$ from our simulation results and validation $E^s$ from $JCMWave$ for the circle and rectangle, respectively.}
\label{error}
 \end{figure}
 
\begin{figure}[htb]
 \centering \subfigure[]{ 
 \label{example1_real_z_neg02} 
  \includegraphics[height=4cm]{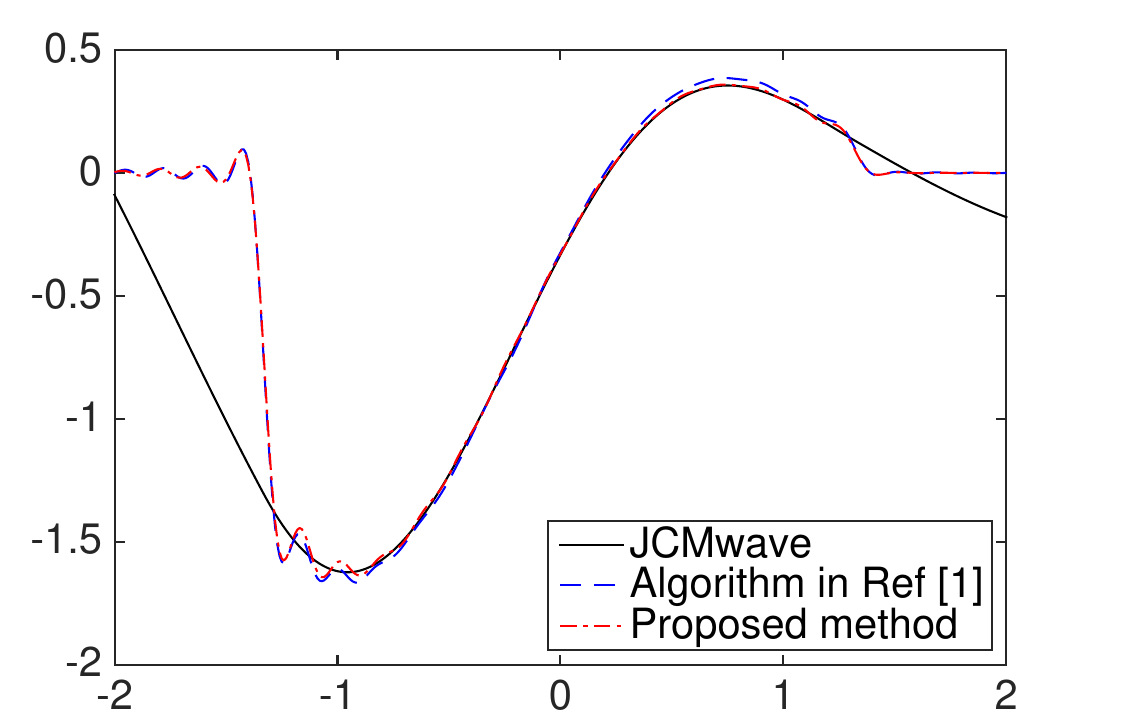}} 
 \subfigure[]{ 
 \label{example1_imag_z_neg02}
 \includegraphics[height=4cm]{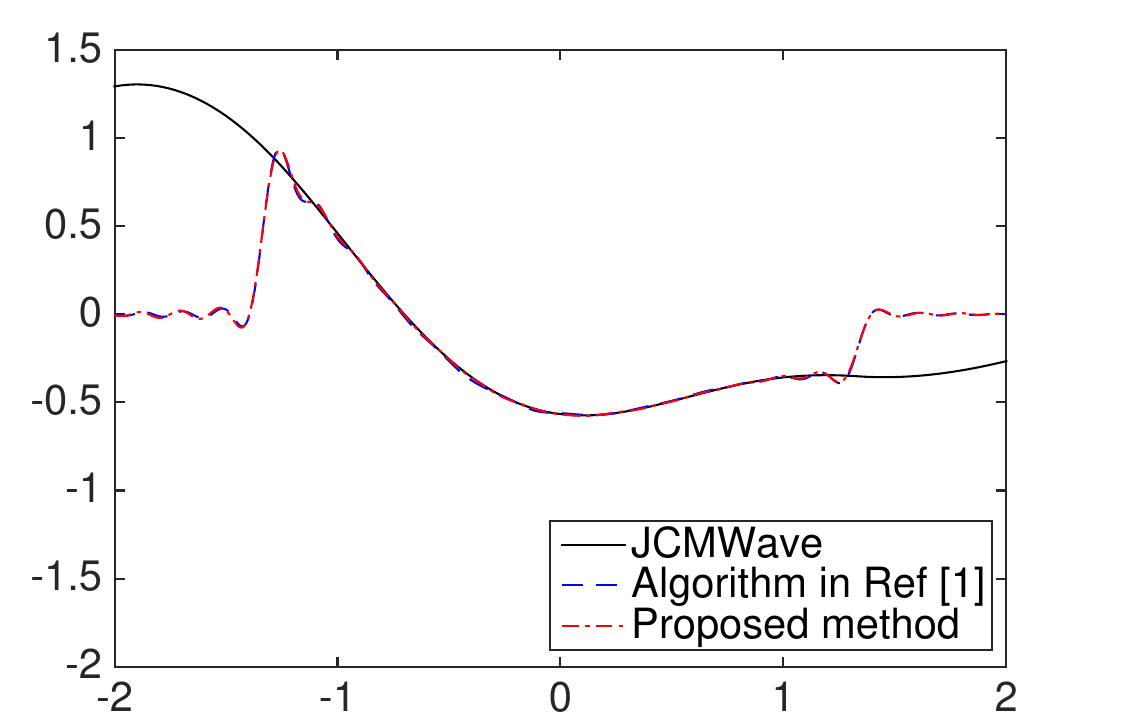}}  
 \caption{The real (shown in (a)) and imaginary (shown in (b)) part of $E^s$ from validation (black solid) and $\chi E^s$ from the algorithm in \cite{Dilz2016The} (blue dashed) and our simulation (red dashdot) at $z = -0.2$ respectively for the circle.}
 \label{example1_z_neg02}
 \end{figure}
 \begin{figure}[htb]
 \centering \subfigure[]{ 
 \label{example2_real_z_0} 
  \includegraphics[height=4cm]{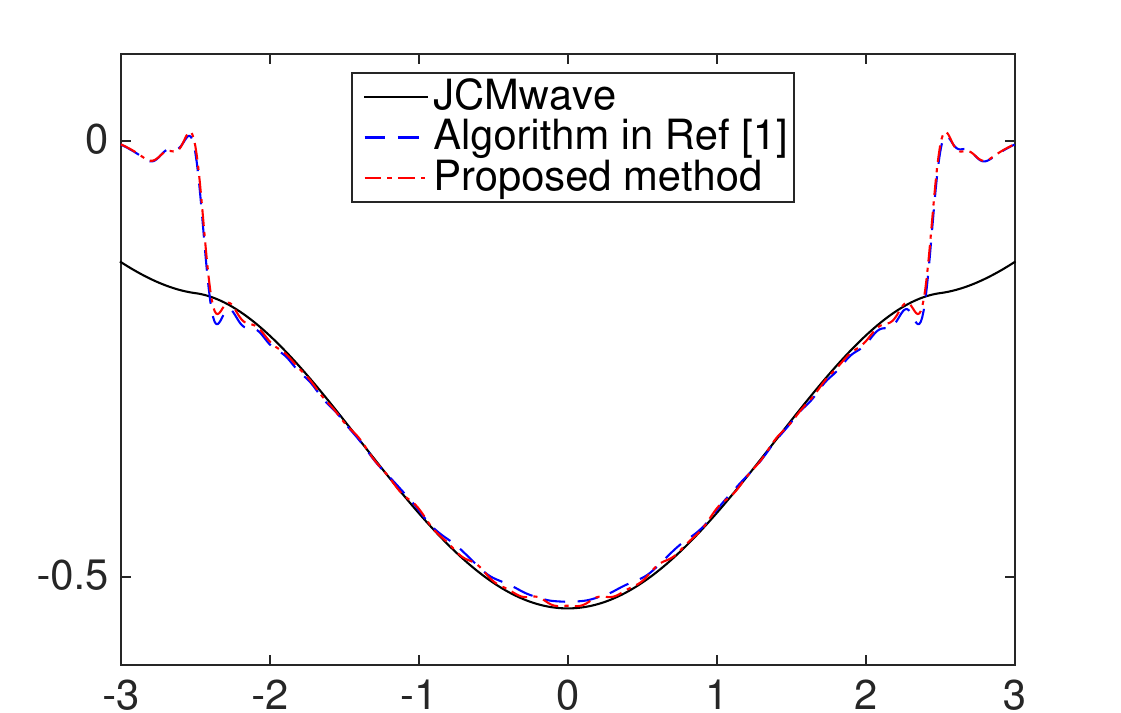}} 
 \subfigure[]{ 
 \label{example2_imag_z_0}
  \includegraphics[height=4cm]{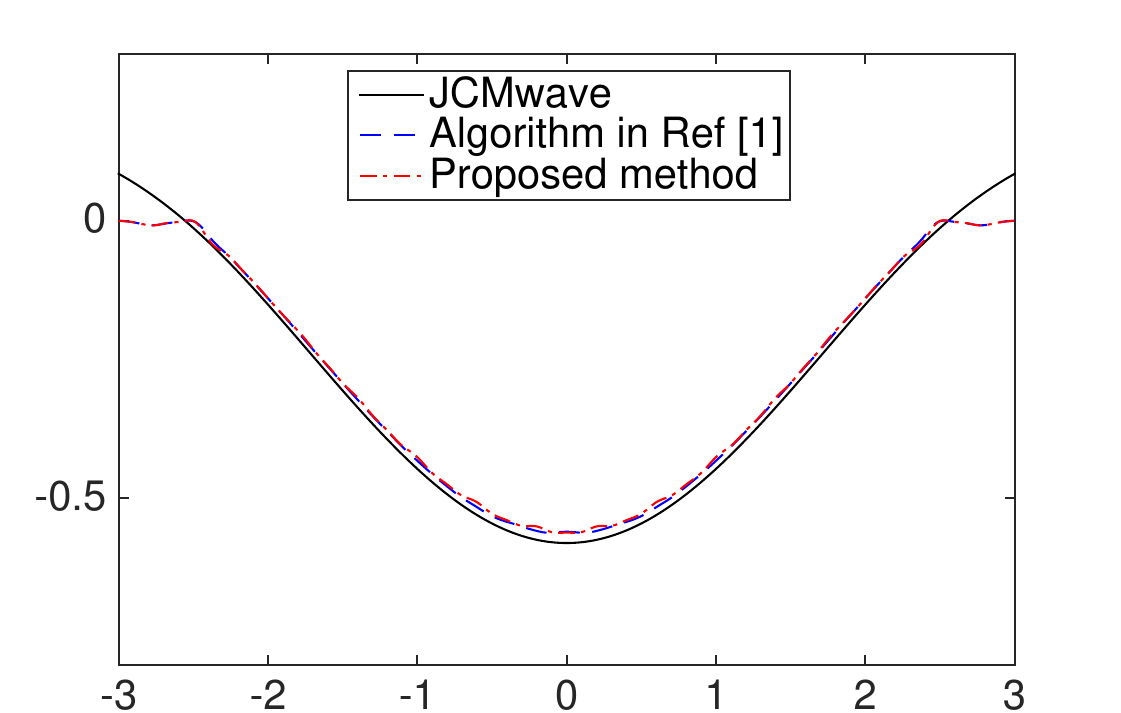}} 
 \caption{The real (shown in (a)) and imaginary (shown in (b)) part of $E^s$ from validation (black solid) and $\chi E^s$ from the algorithm in \cite{Dilz2016The} (blue dashed) and our simulation (red dashdot) at $z = 0$ respectively for the rectangle.}
 \label{example2_z_0}
 \end{figure}

 The real part of the scattered field calculated by our method is shown in Fig.~\ref{realpart}. 
 To validate our results, we utilize the numerical solutions calculated by the JCMWave software package \cite{Burger2013Finite} as the benchmarks. Fig.~\ref{error} presents the difference between the scattered electric field $E^s$ from JCMWave and $\chi E^s$ from our simulations. The large error around the discontinuity is due to the Gibbs phenomenon. Additionally, outside the object, the scattered electric field $E^s$ and the contrast source $\chi E^s$ have a large deviation since the contrast source is zero outside the support of the objects, whereas the scattered electric field is nonzero, which explains the large error outside the objects. It can also be observed in Fig.~\ref{example1_z_neg02} and Fig.~\ref{example2_z_0}. To show the details in the comparison, we plot the real and imaginary part of $E^s$ from JCMWave and $\chi E^s$ from the algorithm described in \cite{Dilz2016The} and our simulation at $z = -0.2$ for the circle in Fig.~\ref{example1_z_neg02} and {at} $z = 0$ for the rectangle in Fig.~\ref{example2_z_0}. 
 These error figures illustrate that our simulations {are} reliable and very {close} to the results obtained by the algorithm in \cite{Dilz2016The}.

\begin{figure}[htb]
 \centering
  \includegraphics[height=2.5cm]{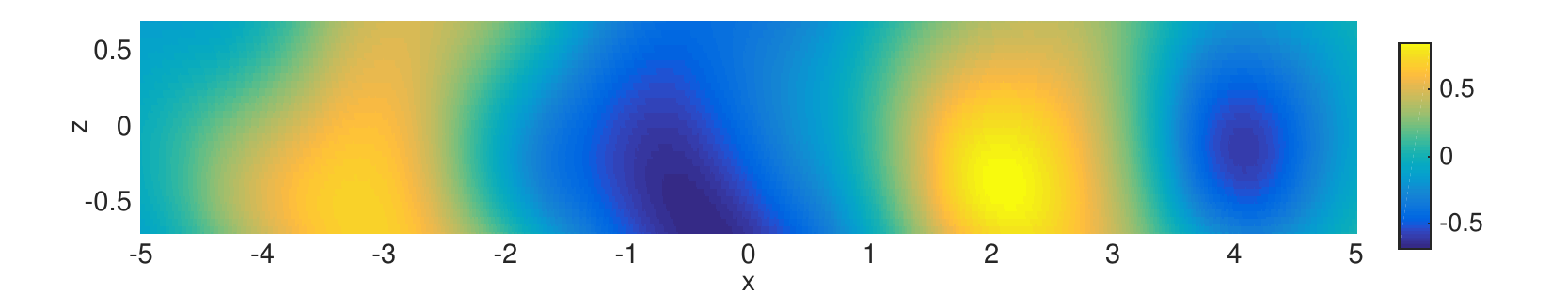}
 \caption{The real parts of $\chi E^s (x,z)$ for the finite grating structure.}
 \label{example3_realpart}
 \end{figure}

\begin{figure}[!h]
 \centering
  \includegraphics[height=2.5cm]{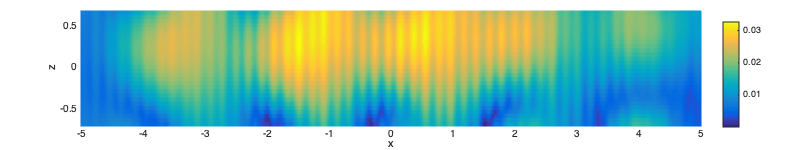}
 \caption{The absolute value of the error between $\chi E^s$ from our simulation results and validation $E^s$ from JCMWave for the finite grating structure.}
 \label{example3_error}
 \end{figure}

\begin{figure}[!h]
 \centering \subfigure[]{ 
 \label{example3_real_z_0} 
 \includegraphics[height=4cm]{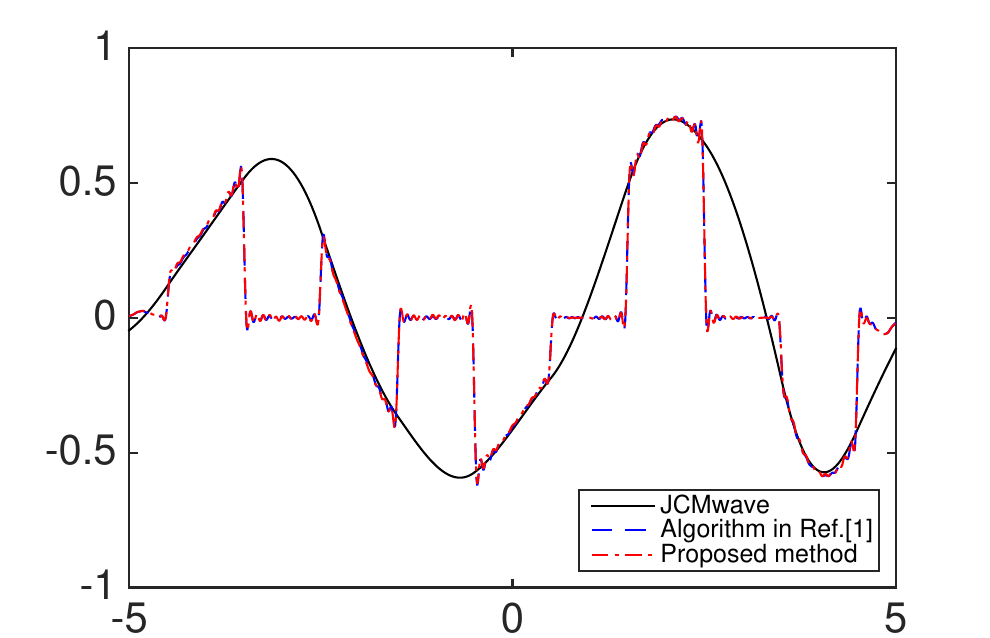}} 
 \subfigure[]{ 
 \label{example3_imag_z_0}
 \includegraphics[height=4cm]{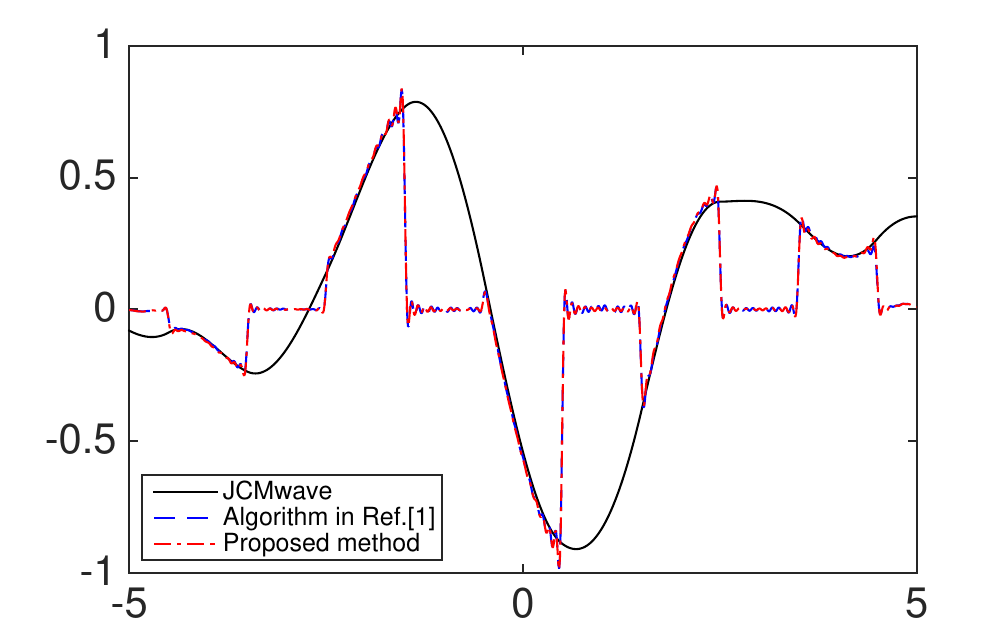}} 
 \caption{The real (shown in (a)) and imaginary (shown in (b)) part of $E^s$ from validation (black solid) and $\chi E^s$ from the algorithm in \cite{Dilz2016The} (blue dashed) and our simulation (red dashdot) at $z = 0$ for the finite grating structure.}
 \label{example3_z_0}
 \end{figure}

 We have also tested our code on the finite grating structure shown in Table \ref{scatters_parameters}. The real part of the scattered field and the error are shown in Fig.~\ref{example3_realpart} and Fig.~\ref{example3_error}. The detailed comparisons of the real and imaginary part of the scattered field at $z=0$ can be seen in Fig.~\ref{example3_real_z_0} and Fig.~\ref{example3_imag_z_0}. 

\subsection{Computational costs}
 To demonstrate the efficiency of our algorithm, we compare the computation times with the algorithm described in \cite{Dilz2016The}. The two algorithms  are both implemented {in} MATLAB R2017b. The parameters of the discretizations for these scatterers are presented in {Table \ref{discretization_parameters}},
 \begin{small}
 \begin{table} [!h]
 \centering
 \begin{tabular}{|c|c|c|c|c|c|c|}
 \hline
 Parameters of the discretization &$X$ & $M$ &$N$  & $\alpha (=\beta)$ & $\Delta$ & $N_k$ \\\hline
  Circle &0.5& 6 & 3 & $\sqrt{2/3}$ & 0.05 & 56 \\\hline
  Rectangle &0.5& 6 & 3 & $\sqrt{2/3}$ & 0.05 & 56\\\hline
   Grating  & 0.5 & 11 & 7 & $\sqrt{2/3}$ & 0.05 & 33  \\\hline
 \end{tabular}
 \caption{The parameters of the discretization for circle, rectangle, and grating, respectively.}
 \label{discretization_parameters}
 \end{table}
 \end{small}

 
 \begin{small}
 \begin{table}[!htb]
  \centering
    \begin{tabular}{|p{2.7cm}|p{1.1cm}|p{1.1cm}|p{1.8cm}|p{1.1cm}|p{1.1cm}|p{1.8cm}|}
    \hline
    \multirow{2}{*}{ \tabincell{l}{Computation\\time (seconds)}  } &
    \multicolumn{3}{c|}{Algorithm Ref \cite{Dilz2016The}} &
    \multicolumn{3}{c|}{ Our proposed method}  \\
    \cline{2-7}
    & MSetT & MSolT  & total time  &   
     MSetT & MSolT  & total time   \\
    \hline
    Circle &532 &0.411 & 533 &265 &0.443  & 267 \\
    \hline
    Rectangle & 536 & 0.353 & 537 & 269  & 0.340 &270 \\
    \hline
     Grating&5755 & 3.138 & 5759 & 1161 & 3.268 & 1165\\
    \hline
  \end{tabular}
  \caption{The computation times of the two algorithms calculating the three scatterers.}
 \label{CPUtimes}
\end{table}
 \end{small}
\noindent where $MSetT$ {indicates} the time for the final coefficient matrix setup, while $MSolT$ means the time for solving the linear equation. 
 
 Table \ref{CPUtimes} {demonstrates} that 
 the ratio of the {computation times} of our simulation and the algorithm in \cite{Dilz2016The} is about $1/2$. As the parameters $M$ and $N$ increase, this ratio drops to $1/4$. {This is mainly because} the inverse matrix calculation of a $(2M+1)(2N+1)(N_k+1)$ matrix when deriving the relation between the Gabor coefficients of the scattered electric field $E^s$ and that of the incoming electric field $E^i$. As a consequence, our proposed method has a greater advantage when simulating larger structures. 

 Next, we briefly discuss the storage requirements in computation. According to Eq.~\eqref{f_spatial_mnl_stk} and Eq.~\eqref{E_spectral_nml_stk}, we utilize two three-dimensional {arrays} to {store} the integrals in {the} spatial and spectral domain. These integrands are continuous and their convergence rates are shown in Eq.~\eqref{limit_spatial} and Eq.~\eqref{limit_spectral}. By observing the integrands, we see that the integral matrices can be just calculated once if the window width $X$ and the oversampling rates $\alpha$, $\beta$ remain unchanged. Then, we can assemble the coefficient matrices of the spatial and spectral part via simple linear computations of the two integral arrays. The scaling of the number of coefficients in these arrays is $O(MNN_k)$.

\section{Conclusions}
\label{conclusions}
 In this paper, we proposed an efficient method to solve the domain integral equation for simulating 2D transverse electric scattering of a finite size object in a homogeneous medium using the Ewald transformation and Gabor frame discretization. The Ewald Green function transformation separates the integrals related to $x$ and $z$, and the Ewald method splits the integral formula of the Green function in two parts, $G_{spatial}$ and $G_{spectral}$. Therefore, the coefficient matrices can be obtained through integral and linear operations without matrix inversion, {which} is the main reason why our proposed method can save computational time. Next we utilize the Fourier {transformation} to deal with the singularity in $G_{spectral}$.  The Gabor frame ensures fast conversion between {the} spatial domain and spectral domain. In addition, we use the summation of modulated Gaussian functions to approximate the dual Gabor window $\eta (x)$ instead of discrete numerical values. As a consequence, we are faced with integrands composed of modulated Gaussian functions and this kind of integrals has an analytic solution. Consequently, the three-dimensional integrals in Eq.~\eqref{f_spatial_mnl_stk} and Eq.~\eqref{E_spectral_nml_stk} {are reduced to one-dimensional integrals.}

Our work in the future is to extend the method to 2D TE scattering in {a} layered medium and the {three-dimensional} case.





\section*{Reference}
\bibliographystyle{model1-num-names}
\bibliography{sample}

\begin{thebibliography}{23}
\expandafter\ifx\csname natexlab\endcsname\relax\def\natexlab#1{#1}\fi
\providecommand{\bibinfo}[2]{#2}
\ifx\xfnm\relax \def\xfnm[#1]{\unskip,\space#1}\fi
\bibitem[{Zwamborn and van~den Berg(1991)}]{PP1991Aweak}
\bibinfo{author}{P.~Zwamborn}, \bibinfo{author}{P.~M. van~den Berg},
\newblock \bibinfo{title}{A weak form of the conjugate gradient {FFT} method
  for two-dimensional {TE} scattering problems},
\newblock \bibinfo{journal}{IEEE Transactions on Microwave Theory and
  Techniques} \bibinfo{volume}{39} (\bibinfo{year}{1991})
  \bibinfo{pages}{953--960}.
\bibitem[{Bleszynski et~al.(1996)Bleszynski, Bleszynski, and
  Jaroszewicz}]{AIM1996EB}
\bibinfo{author}{E.~Bleszynski}, \bibinfo{author}{M.~Bleszynski},
  \bibinfo{author}{T.~Jaroszewicz},
\newblock \bibinfo{title}{{AIM}: {A}daptive integral method for solving
  large-scale electromagnetic scattering and radiation problems},
\newblock \bibinfo{journal}{Radio Science} \bibinfo{volume}{31}
  (\bibinfo{year}{1996}) \bibinfo{pages}{1225--1251}.
\bibitem[{Lalanne et~al.(2007)Lalanne, Besbes, Hugonin, Haver, Janssen,
  Nugrowati, Xu, Pereira, Urbach, and Nes}]{Lalanne2007Numerical}
\bibinfo{author}{P.~Lalanne}, \bibinfo{author}{M.~Besbes},
  \bibinfo{author}{J.~P. Hugonin}, \bibinfo{author}{S.~V. Haver},
  \bibinfo{author}{O.~T.~A. Janssen}, \bibinfo{author}{A.~M. Nugrowati},
  \bibinfo{author}{M.~Xu}, \bibinfo{author}{S.~F. Pereira},
  \bibinfo{author}{H.~P. Urbach}, \bibinfo{author}{A.~S. V.~D. Nes},
\newblock \bibinfo{title}{Numerical analysis of a slit-groove diffraction
  problem},
\newblock \bibinfo{journal}{Journal of the European Optical Society Rapid
  Publications} \bibinfo{volume}{2} (\bibinfo{year}{2007}).
\bibitem[{Zwamborn and van~den Berg(1992)}]{Zwamborn1992The}
\bibinfo{author}{P.~Zwamborn}, \bibinfo{author}{P.~M. van~den Berg},
\newblock \bibinfo{title}{The three dimensional weak form of the conjugate
  gradient fft method for solving scattering problems},
\newblock \bibinfo{journal}{IEEE Transactions on Microwave Theory and
  Techniques} \bibinfo{volume}{40} (\bibinfo{year}{1992})
  \bibinfo{pages}{1757--1766}.
\bibitem[{Phillips and White(1994)}]{Phillips1994Proceedings}
\bibinfo{author}{J.~R. Phillips}, \bibinfo{author}{J.~K. White},
\newblock \bibinfo{title}{Efficient capacitance extraction of {3D} structures
  using generalized pre-corrected fft methods},
\newblock in: \bibinfo{booktitle}{Proceedings of 1994 IEEE Electrical
  Performance of Electronic Packaging}, pp. \bibinfo{pages}{253--256}.
\bibitem[{Dilz and {van Beurden}(2016)}]{Dilz2016The}
\bibinfo{author}{R.~J. Dilz}, \bibinfo{author}{M.~C. {van Beurden}},
\newblock \bibinfo{title}{The {G}abor frame as a discretization for the 2{D}
  transverse-electric scattering-problem domain integral equation},
\newblock \bibinfo{journal}{Progress in Electromagnetics Research B}
  \bibinfo{volume}{69} (\bibinfo{year}{2016}) \bibinfo{pages}{117--136}.
\bibitem[{Bastiaans(1981)}]{Bastiaans1981A}
\bibinfo{author}{M.~J. Bastiaans},
\newblock \bibinfo{title}{A sampling theorem for the complex spectrogram, and
  {G}abor's expansion of a signal in {G}aussian elementary signals},
\newblock \bibinfo{journal}{Optical Engineering} \bibinfo{volume}{20}
  (\bibinfo{year}{1981}) \bibinfo{pages}{594--598}.
\bibitem[{Bastiaans(1995)}]{Bastiaans1995Gabor}
\bibinfo{author}{M.~J. Bastiaans},
\newblock \bibinfo{title}{{G}abor's expansion and the {Z}ak transform for
  continuous-time and discrete-time signals: Critical sampling and rational
  oversampling}  (\bibinfo{year}{1995}).
\bibitem[{Feichtinger and Strohmer(1998)}]{Feichtinger1998Gabor}
\bibinfo{author}{H.~G. Feichtinger}, \bibinfo{author}{T.~Strohmer},
  \bibinfo{title}{Gabor Analysis and Algorithms: Theory and Applications},
  \bibinfo{publisher}{Birkhauser Boston}, \bibinfo{year}{1998}.
\bibitem[{Werther et~al.(2005)Werther, Eldar, and Subbanna}]{Werther2005Dual}
\bibinfo{author}{T.~Werther}, \bibinfo{author}{Y.~C. Eldar},
  \bibinfo{author}{N.~K. Subbanna},
\newblock \bibinfo{title}{Dual {G}abor frames: theory and computational
  aspects},
\newblock \bibinfo{journal}{IEEE Transactions on Signal Processing}
  \bibinfo{volume}{53} (\bibinfo{year}{2005}) \bibinfo{pages}{4147--4159}.
\bibitem[{Janssen(1994)}]{Janssen1994Signal}
\bibinfo{author}{A.~J. E.~M. Janssen},
\newblock \bibinfo{title}{Signal analytic proofs of two basic results on
  lattice expansions},
\newblock \bibinfo{journal}{Applied \& Computational Harmonic Analysis}
  \bibinfo{volume}{1} (\bibinfo{year}{1994}) \bibinfo{pages}{350--354}.
\bibitem[{Ewald(1921)}]{Ewald1921Die}
\bibinfo{author}{P.~P. Ewald},
\newblock \bibinfo{title}{Die berechnung optischer und elektrostatischer
  gitterpotentiale},
\newblock \bibinfo{journal}{Annalen Der Physik} \bibinfo{volume}{369}
  (\bibinfo{year}{1921}) \bibinfo{pages}{253--287}.
\bibitem[{Ewald(1970)}]{Ewald1970On}
\bibinfo{author}{P.~P. Ewald},
\newblock \bibinfo{title}{On the foundations of crystal optics. part 1.
  dispersion theory. part 2. theory of reflection and refraction}
  (\bibinfo{year}{1970}).
\bibitem[{Jordan et~al.(1986)Jordan, Richter, and Sheng}]{Jordan1986An}
\bibinfo{author}{K.~E. Jordan}, \bibinfo{author}{G.~R. Richter},
  \bibinfo{author}{P.~Sheng},
\newblock \bibinfo{title}{An efficient numerical evaluation of the {G}reen's
  function for the {H}elmholtz operator on periodic structures},
\newblock \bibinfo{journal}{Journal of Computational Physics}
  \bibinfo{volume}{63} (\bibinfo{year}{1986}) \bibinfo{pages}{222--235}.
\bibitem[{Mathis and Peterson(1996)}]{Mathis1996A}
\bibinfo{author}{A.~W. Mathis}, \bibinfo{author}{A.~F. Peterson},
\newblock \bibinfo{title}{A comparison of acceleration procedures for the
  two-dimensional periodic {G}reen's function},
\newblock \bibinfo{journal}{IEEE Transactions on Antennas \& Propagation}
  \bibinfo{volume}{44} (\bibinfo{year}{1996}) \bibinfo{pages}{567--571}.
\bibitem[{Capolino et~al.(2005)Capolino, Wilton, and
  Johnson}]{Capolino2005Efficient}
\bibinfo{author}{F.~Capolino}, \bibinfo{author}{D.~R. Wilton},
  \bibinfo{author}{W.~A. Johnson},
\newblock \bibinfo{title}{Efficient computation of the 2{D} {G}reen's function
  for 1{D} periodic structures using the {E}wald method},
\newblock \bibinfo{journal}{IEEE Transactions on Antennas \& Propagation}
  \bibinfo{volume}{53} (\bibinfo{year}{2005}) \bibinfo{pages}{2977--2984}.
\bibitem[{Capolino et~al.(2007)Capolino, Wilton, and
  Johnson}]{Capolino2007Efficient}
\bibinfo{author}{F.~Capolino}, \bibinfo{author}{D.~R. Wilton},
  \bibinfo{author}{W.~A. Johnson},
\newblock \bibinfo{title}{Efficient computation of the 3{D} {G}reen's function
  for the {H}elmholtz operator for a linear array of point sources using the
  {E}wald method},
\newblock \bibinfo{journal}{Journal of Computational Physics}
  \bibinfo{volume}{223} (\bibinfo{year}{2007}) \bibinfo{pages}{250--261}.
\bibitem[{Komanduri et~al.(2016)Komanduri, Jackson, Capolino, and
  Wilton}]{Komanduri20161}
\bibinfo{author}{V.~R. Komanduri}, \bibinfo{author}{D.~R. Jackson},
  \bibinfo{author}{F.~Capolino}, \bibinfo{author}{D.~R. Wilton},
\newblock \bibinfo{title}{1{D} periodic {G}reen's function for leaky and
  complex waves using the {E}wald method},
\newblock \bibinfo{journal}{IEEE Transactions on Antennas \& Propagation}
  \bibinfo{volume}{64} (\bibinfo{year}{2016}) \bibinfo{pages}{4703--4712}.
\bibitem[{Dilz(2017)}]{RJDilz2017spatialspectral}
\bibinfo{author}{R.~Dilz}, \bibinfo{title}{A spatial spectral domain integral
  equation solver for electromagnetic scattering in dielectric layered media},
  Ph.D. thesis, Eindhoven University of Technology, \bibinfo{year}{2017}.
\bibitem[{Arens et~al.(2013)Arens, Kai, Schmitt, and
  Lechleiter}]{Arens2013Analysing}
\bibinfo{author}{T.~Arens}, \bibinfo{author}{S.~Kai},
  \bibinfo{author}{S.~Schmitt}, \bibinfo{author}{A.~Lechleiter},
\newblock \bibinfo{title}{Analysing {E}wald's method for the evaluation of
  {G}reen's functions for periodic media},
\newblock \bibinfo{journal}{{IMA} Journal of Applied Mathematics}
  \bibinfo{volume}{78} (\bibinfo{year}{2013}) \bibinfo{pages}{405--431}.
\bibitem[{Gradshteyn and Ryzhik(2014)}]{2014Table}
\bibinfo{author}{I.~S. Gradshteyn}, \bibinfo{author}{I.~M. Ryzhik},
  \bibinfo{title}{Table of Integrals, Series, and Products (Eighth Edition)},
  \bibinfo{publisher}{Academic press}, \bibinfo{year}{2014}.
\bibitem[{Abramowitz et~al.(1972)Abramowitz, Stegun, and
  Romain}]{M1972Handbook}
\bibinfo{author}{M.~Abramowitz}, \bibinfo{author}{I.~A. Stegun},
  \bibinfo{author}{J.~E. Romain}, \bibinfo{title}{Handbook of Mathematical
  Functions with Formulas, Graphs, and Mathematical Tables},
  \bibinfo{publisher}{New York: Dover}, \bibinfo{year}{1972}.
\bibitem[{Burger et~al.(2013)Burger, Lin, Pomplun, and
  Schmidt}]{Burger2013Finite}
\bibinfo{author}{S.~Burger}, \bibinfo{author}{Z.~Lin},
  \bibinfo{author}{J.~Pomplun}, \bibinfo{author}{F.~Schmidt},
\newblock \bibinfo{title}{Finite-element based electromagnetic field
  simulations: benchmark results for isolated structures},
\newblock in: \bibinfo{booktitle}{SPIE Photomask Technology}.

\end{thebibliography}







\nomenclature{$X,K,\alpha,\beta$}{Gabor frame parameters (constants)}
\nomenclature{$x,z,x',z'$}{Coordinate variables}
\nomenclature{$k_0$}{Wavenumber in vacuum}
\nomenclature{$k_x$}{Wavenumber in $x$-direction}
\nomenclature{$R$}{The distance between observation point $(x,z)$ and the source element $(x',z')$}
\nomenclature{$w$}{Continuous real variable}
\nomenclature{$\xi,\zeta$}{Complex variables}
\nomenclature{$\phi,f,g,h$}{functions in the spatial domain}
\nomenclature{$\widetilde{f},\widetilde{g},\widetilde{h}$}{functions in the spectral domain}
\nomenclature{$\begin{array}{cc} m,n,l,s,t,k, \\ p,q,u,v \end{array}$}{Integer variables as the index of Gabor coefficients}
\nomenclature{$N_\theta $}{The upper limit of absolute value of $\theta$, where $\theta$ can be any one in  $\{m,n,l,s,t,k,p,q,u,v\}$}

\end{document}